\documentclass[10pt]{article}                                                                                                                           
\usepackage{amssymb}                                                            
\usepackage{amsmath,amsthm} 
\usepackage{mathrsfs}                                                            
\usepackage[all]{xy}
\usepackage[dvips]{graphicx}
\usepackage{wrapfig} 

\newcommand{\ud}{\mathrm{d}}
\newcommand{\cl}{{\cal L}}
\newcommand{\p}{\partial}
\newtheorem{prop}{Proposition}[section]
\newtheorem{col}{Corollary}[section]

\newtheorem{theorem}{Theorem}[section]

\numberwithin{equation}{section}

\title{\bf{Quantum Group $GL_q(2)$ and Quantum Laplace Operator via Semi-infinite Cohomology\\
}} 
\author{\vspace{5mm}  
Igor B. Frenkel\footnote{igor.frenkel@yale.edu}
\hspace{5mm}Anton M. Zeitlin\footnote{anton.zeitlin@yale.edu, http://math.yale.edu/$\sim$az84 http://www.ipme.ru/zam.html}\\
Department of Mathematics,\\
Yale University,\\
442 Dunham Lab, 10 Hillhouse Avenue,\\
New Haven, CT 06511}

\begin{document}
\maketitle
\begin{abstract}
We construct the quantum group $GL_q(2)$ as the semi-infinite cohomology of the tensor product of two braided vertex operator algebras based on the algebra $W_2$ with complementary central charges $c+\bar{c}=28$. The conformal field theory version of the Laplace operator on the quantum group is also obtained.
\end{abstract}

\section{Introduction}
We have shown in \cite{fz} that the quantum group $SL_q(2)$ admits a realization as the semi-infinite cohomology of the Virasoro algebra with coefficients in the tensor product of two braided vertex operator algebras (VOA)  with complementary central charges $c+\bar{c}=26$. At the end of the paper we have discussed extensions of our results to higher ranks and more general $W$-algebras (see e.g. \cite{FFr} and references therein). In the present paper we consider a simplest generalization of this type, namely the algebra $W_2$, which is the semi-direct sum of the Virasoro and the Heisenberg algebra. We prove that the semi-infinite cohomology of $W_2$ in the tensor product of two appropriate braided vertex operator algebras yields the quantum group $GL_q(2)$. The central charges of the two braided  VOAs add up to the critical value for $W_2$, namely $c+\bar{c}=28$.

It is well known that one can develop the calculus of differential forms on $GL_q(2)$ and define the $q$-analogues of partial derivatives and the quantum version of Laplace operator. Since $GL_q(2)$ is realized as semi-infinite cohomology of $W_2$, one can look for a "lift" of various structures on the quantum group to the corresponding braided VOA. In this paper we do find a simple operator on the braided VOA that commutes with the action of 
$W_2$  and induces the quantum Laplace operator on the cohomology. This suggests that other structures of noncommutative geometry might admit equally simple and natural realization in terms of chiral two-dimensional conformal field theory. 

The structure of the paper is as follows. In Section 2 we recall the basic facts about representation theory of 
$U_q(gl(2))$.  In Section 3 we discuss $W_2$ and related semi-infinite cohomology complex. In Section 4 we construct the intertwining operators for $W_2$ and related braided VOA (see \cite{styrkas}, \cite{fz}). 
In Section 5 we find realization of $GL_q(2)$ as the semi-infinite cohomology of certain braided VOA and write a simple formula for the quantum Laplace operator using a Fock space realization of $W_2$ modules. 

\section{$\mathbf{U_q(gl(2))}$, its representations and intertwining operators}
Let $U_q(gl(2))$ be the Hopf algebra over $\mathbb{C}(q)$ with generators $E, F, q^{\pm H}, I$
and commutation relations:
\begin{eqnarray}
 q^HE&=&q^2Eq^H,\nonumber\\
q^HF&=&q^{-2}Fq^H,\nonumber \\
{[}E,F] &=&\frac{q^H-q^{-H}}{q-q^{-1}},\nonumber
\end{eqnarray}
and $I$ is a central element. The comultiplication is given by
\begin{eqnarray}
\Delta(I)&=&I\otimes 1+1\otimes I,\nonumber\\
\Delta(q^H)&=&q^{H}\otimes q^{H}\nonumber\\
\Delta(E)&=&E\otimes q^H+1\otimes E,\nonumber\\
\Delta(F)&=&F\otimes 1 + q^{-H}\otimes F.
\end{eqnarray}
The universal R-matrix for $U_q(gl(2))$, which is an element of a certain completion of 
$U_q(gl(2))\otimes U_q(gl(2))$, is given by:
\begin{eqnarray}
R&=&C\Theta, \qquad C=q^{{I\otimes I}}q^{\frac{H\otimes H}{2}},\nonumber\\
\Theta&=&\sum_{k\geqslant 0}q{^{k(k-1)/2}\frac{(q-q^{-1})^k}{[k]!}E^k\otimes F^k},
\end{eqnarray}
where $[n]=\frac{q^n-q^{-n}}{q-q^{-1}}$ and $[n]!=[1][2]\ldots[n]$.

For any given pair $V,W$ of representations, R-matrix gives the following commutativily isomorphism:
$\check{R}=PR: V\otimes W \rightarrow W\otimes V$, where $P$ is a permutation: $P(v\otimes w)=w\otimes v$.

Let $\lambda\in \mathbb{Z}_+$ and $k\in \mathbb{Z}$.
We denote by $V_{\lambda,k}$ the finite dimensional irreducible representation of $U_q(gl(2))$, such that $V_{\lambda,k}$ is the highest weight representation with the highest weight $\lambda$ for the subalgebra $U_q(sl(2))$ of $U_q(gl(2))$, and 
$k$ is the eigenvalue of central element $I$. 

One can construct the intertwining operators for finite dimensional representations, i.e. elements of $Hom (V_{\lambda , l}\otimes V_{\mu ,m}, V_{\nu ,n})$ and 
$Hom (V_{\nu ,n}, V_{\lambda , l}\otimes V_{\mu ,m})$. 
The following Proposition holds:
\begin{prop}
Let $\lambda, \mu, \nu\in\mathbb{Z}_+$ and $l,m,n\in \mathbb{Z}$. Then 
$dim Hom (V_{\lambda, l}\otimes V_{\mu,m}, V_{\nu,n})=1$ iff $\lambda+\mu\ge\nu\ge|\lambda-\mu|$ and $n=m+l$. Otherwise $dim Hom (V_{\lambda, l}\otimes V_{\mu,m}, V_{\nu,n})=0$.
\end{prop}
\noindent {\bf Proof.} Proof follows from a similar fact from $U_q(sl(2))$ representation theory.

\hfill$\blacksquare$\bigskip

\noindent The same statement holds for intertwiners from 
$Hom (V_{\nu ,n}, V_{\lambda , l}\otimes V_{\mu ,m})$.
\noindent The next proposition gives quadratic relations between intertwiners and will be very crucial in the following.
\begin{prop}
Let $\lambda_i\in \mathbb{Z}_+$, $l_i\in \mathbb{Z}$ $(i=0,1,2,3)$ . 
Then there exists an invertible operator 
\begin{displaymath}
B\left[\begin{array}{cc}
\lambda_0, l_0 ; \lambda_1, l_1\\
\lambda_2, l_2 ; \lambda_3,l_3
\end{array}
\right]
\end{displaymath}
such that the following diagram is commutative:

\begin{displaymath}
\xymatrixcolsep{50pt}
\xymatrixrowsep{10pt}
\xymatrix{
{\begin{array}{l}
\oplus_\rho \big(Hom(V_{\rho, r},V_{\lambda_1, l_1}\otimes V_{\lambda_2,l_2})\\
\otimes Hom(V_{\lambda_0, l_0}, V_{\rho, r}\otimes  V_{\lambda_3, l_3})\big)
\end{array}
\ar[d]_{i}}
\ar[r]^{
B\left[\begin{array}{cc}
\lambda_0,l_0; \lambda_1,l_1\\
\lambda_2,l_2 ;\lambda_3,l_3
\end{array}
\right]
}&  
{\begin{array}{r}
\oplus_{\xi, k} \big(Hom(V_{\xi,k},V_{\lambda_1, l_1}\otimes V_{\lambda_3, l_3})\\
\otimes Hom(V_{\lambda_0, l_0}, V_{\xi, k}\otimes V_{\lambda_2, l_2})\big)
\end{array}
\ar[d]_{i}
}
\\
Hom(V_{\lambda_0, l_0}, V_{\lambda_1, l_1}\otimes V_{\lambda_2, l_2}\otimes V_{\lambda_3, l_3})
\ar[r]^{PR}& 
Hom(V_{\lambda_0, l_0}, V_{\lambda_1, l_1}\otimes V_{\lambda_3, l_3}\otimes V_{\lambda_2, l_2}), 
}
\end{displaymath}

where $|\lambda_1+\lambda_2|\ge \rho\ge |\lambda_1-\lambda_2|$, $|\lambda_3+\rho|\ge \lambda_0\ge |\lambda_3-\rho|$, $|\lambda_1+\lambda_3|\ge \xi\ge |\lambda_1-\lambda_3|$, $|\lambda_2+\xi|\ge \lambda_0\ge |\lambda_2-\xi|$, and $i$ is an isomorphism.
\end{prop}

\noindent This fact follows from a similar fact from the representation theory of $U_q(sl(2))$. Moreover, when the braiding matrix 
 \begin{eqnarray}
 B\left[\begin{array}{cc}
\lambda_0,l_0; \lambda_1,l_1\\
\lambda_2,l_2 ;\lambda_3,l_3
\end{array}
\right]
\end{eqnarray} 
 is not equal to zero, it is equal to  
\begin{eqnarray}
 q^{l_2l_3}B^V\left[\begin{array}{cc}
\lambda_0; \lambda_1\\
\lambda_2 ;\lambda_3
\end{array}
\right],
\end{eqnarray} 
where $B^V$ is the braiding matrix of $U_q(sl(2))$. Let's now denote  by 
$\phi^{\nu,n}_{\lambda,l;\mu,m}$ the generating element of 
$Hom (V_{\lambda , l}\otimes V_{\mu ,m}, V_{\nu ,n})$ and by $\phi_{\nu,n}^{\lambda,l;\mu,m}$ the generating element of 
$Hom (V_{\nu ,n}, V_{\lambda , l}\otimes V_{\mu ,m})$. Then Proposition 2.2. gives the following quadratic relations between intertwining operators:
\begin{eqnarray}\label{intfd}
&&(1\otimes PR)\phi_\rho^{\lambda_1,l_1;\lambda_2,l_2}\phi_{\lambda_0,l_0}
^{\rho,r;\lambda_3,l_3}=\nonumber\\
&&\sum_{\xi}B_{\rho,r;\xi,k}\left[\begin{array}{cc}
\lambda_0,l_0 & \lambda_1,l_1\\
\lambda_2,l_2 & \lambda_3,l_3
\end{array}
\right]
\phi_{\xi,k}^{\lambda_1,l_1;\lambda_3,l_3}\phi_{\lambda_0,l_0}^{\xi,k;\lambda_2,l_2},\\
&&\label{intfdnew}
\phi_{\rho,r;\lambda_3,l_3}^{\lambda_0,l_0}\phi_{\lambda_1,l_1\lambda_2,l_2}^\rho(1\otimes PR)=\nonumber\\
&&\sum_{\xi,k}B_{\xi,k;\rho,r}\left[\begin{array}{cc}
\lambda_0,l_0 & \lambda_1,l_1\\
\lambda_2,l_2 & \lambda_3,l_3
\end{array}
\right]
\phi_{\xi,k;\lambda_2,l_2}^{\lambda_0,l_0}\phi_{\lambda_1l_1,\lambda_3,l_3}^{\xi,k}.
\end{eqnarray}

\section{Virasoro, $W_2$ algebras and semi-infinite cohomology: basic facts}
{\bf 3.1. Virasoro algebra and Feigin-Fuks representation.}
The Virasoro algebra
\begin{equation}
[L_n,L_m]=(n-m)L_{m+n}+\frac{\hat c}{12}(n^3-n)\delta_{n,-m}
\end{equation}
has been extensively studied for many years.
Here we need only basic facts. Let us denote by $M_{\hat c,h}$ and $V_{\hat c,h}$ the Verma module and irreducible module
(with highest weight $h$), correspondingly. 
Throughout the paper we will consider only generic values of $\hat c$.  This means that the central charge $\hat{c}$ is parametrized in the following way:
 \begin{equation}
\hat c=13-6(\varkappa+\frac{1}{\varkappa}),
\end{equation}
where parameter $\varkappa\in \mathbb{R}\backslash \mathbb{Q}$. 
Then we have the following proposition (see e.g. \cite{fb} and references therein).
\begin{prop}
For generic value of $\hat c$, Verma module $M_{\hat c,h}$ has a unique singular vector in the case
if $h=h_{m,n}$, where
\begin{equation}\label{virsing}
h_{m,n}=\frac{1}{4}(m^2-1)\varkappa+\frac{1}{4}(n^2-1)\varkappa^{-1}-\frac{1}{2}(mn-1).
\end{equation} 
This singular vector occurs on the level mn, i.e. the value of 
$L_0$ is $h_{m,n}+mn$.
\end{prop}
In the following we will be interested in the modules with
$h=h_{1,n}=\Delta(\lambda)$, where $\lambda=n-1$, $\Delta(\lambda)=-\frac{\lambda}{2}+\frac{\lambda(\lambda+2)}{4\varkappa}$.

\begin{col}Let $\hat c$ be generic and $\lambda \ge 0$, then
$V_{\Delta(\lambda),\hat c}=M_{\Delta(\lambda),\hat c}/M_{\Delta(\lambda)+\lambda+1,\hat c}$, where
$V_{\Delta(\lambda),\hat c}$ is the irreducible Virasoro module with the highest
weight $\Delta(\lambda)$. For $\lambda < 0$ and generic values of $c$ the irreducible module is isomorphic to Verma one, namely, $V_{\Delta(\lambda),\hat c}=M_{\Delta(\lambda),\hat c}$.
\end{col}
Let us consider the Heisenberg algebra
\begin{equation}
[\alpha_n,\alpha_m]=2\varkappa m\delta_{n+m,0}
\end{equation}
and denote by $F_{\lambda,\varkappa}$ the Fock module associated to this algebra. Namely, 
$F_{\lambda,\varkappa}=S(\alpha_{-1},\alpha_{-2},\dots)\otimes\mathbf{1}_\lambda$, such that
$a_n\mathbf{1}_\lambda=0$ if $n>0$ and $a_0\mathbf{1}_\lambda=\lambda\mathbf{1}_\lambda (\lambda\in\mathbb{C})$.
It is well known (see e.g. \cite{fb}, \cite{FHL} and references therein) that $F_{0,\varkappa}$ gives rise to the vertex operator algebra, generated by the field
$\alpha(z)=\sum a_nz^{-n-1}$, such that $\deg(a(z))=1$ which has the following operator product expansion (OPE):
\begin{equation}
\alpha(z)\alpha(w)\sim \frac{2\varkappa}{(z-w)^2}.
\end{equation}
We will denote this vertex algebra as $F_{0,\varkappa}(\alpha)$. The following Proposition holds.
\begin{prop}\label{feiginfuks}
Vertex algebra $F_{0,\varkappa}(a)$ has a vertex operator algebra structure, where the vertex operator, corresponding to the Virasoro 
element, is given by the following formula:
\begin{equation}
L(z)=\frac{1}{4\varkappa}:\alpha(z)^2:+\frac{\varkappa-1}{2\varkappa}\p \alpha(z),
\end{equation} 
such that $L(z)=\sum_nL_nz^{-n-2}$ and $L_n$ satisfy Virasoro algebra
relations with the central charge $\hat c=13-6(\varkappa+\frac{1}{\varkappa})$.
\end{prop}

\noindent {\bf 3.2. $W_2$ as an extension of Virasoro algebra.}
$W_2$ is the semi-direct sum of Virasoro and Heisenberg algebras. The commutation relations are:
\begin{eqnarray}\label{w2}
&&[{\cal L}_n,{\cal L}_m]=(n-m)\cl_{m+n}+\frac{c}{12}(n^3-n)\delta_{n,-m},\nonumber\\
&&[a_n,a_m]=2\eta m\delta_{n+m,0},\nonumber\\
&&[\cl_n,a_m]=-ma_{n+m}.
\end{eqnarray}
The pair of central charges $(c,\eta)$ determine the algebra. Throughout this paper we will use the following parameterization of the central charge:
 \begin{eqnarray}
 c=14-6(\varkappa+\frac{1}{\varkappa}).
 \end{eqnarray}
 Moreover, we will require, that $\varkappa\in \mathbb{R}\backslash \mathbb{Q}$.   
However, for any representation of $W_2$ algebra one can find a 
$direct$ sum of Virasoro and Heisenberg algebras acting in the same representation. Let us define generators 
\begin{eqnarray}
L_n\equiv\cl_n-\frac{1}{2\eta}\sum^{\infty}_{m=-\infty} :a_{n-m}a_m:,
\end{eqnarray}
where symbol $::$ stands for standard Fock space normal ordering. 
Then the following Proposition holds. 
\begin{prop}Let $\cl_n, a_m$ generate algebra $W_2$ with central charges $(c,\eta)$, then $L_n$, $a_n$ satisfy the following commutation relations:
\begin{eqnarray}
&&[L_n,L_m]=(n-m)L_{m+n}+\frac{\hat c}{12}(n^3-n)\delta_{n,-m},\nonumber\\
&&[a_n,a_m]=2\eta m\delta_{n+m,0},\nonumber\\
&&[L_n,a_m]=0,
\end{eqnarray}
where $\hat{c}=c-1$.
\end{prop}
We will be interested in the highest weight modules 
for $W_2$ algebra 
\begin{eqnarray}
W_{\Delta(\lambda),k}^{\varkappa,\eta}=V_{\Delta(\lambda),\varkappa}\otimes F_{k,\eta}
\end{eqnarray}
which are the tensor products of the irreducible highest weight representation of Virasoro algebra $V_{\Delta(\lambda),\varkappa}$, generated by $L_n$ and Fock module $F_{k,\eta}$ for the Heisenberg algebra, generated by $a_n$.\\

\noindent {\bf 3.3. Ghost VOAs and semi-infinite cohomology for $W_2$.} 
In this section we will show how to reduce the semi-infinite cohomology of $W_2$ to the semi-infinite cohomology 
of the Virasoro algebra \cite{fgz}. 

In the case of $W_2$, semi-infinite forms can be realized by means of the following Heisenberg superalgebras:
\begin{eqnarray}
&&\{\psi_n,\chi_m\}=\delta_{n+m,0}, \nonumber\\
&&\{b_n, c_m\}=\delta_{n+m,0}, \quad n,m \in \mathbb{Z}.
\end{eqnarray}
One can construct Fock modules $\Lambda$, $\Lambda'$ in the following way:
\begin{eqnarray}
\Lambda&=&\mathbb{C}\{b_{-n_1}\dots b_{-n_k}c_{-m_1}\dots c_{-m_\ell} \mathbf{1};\nonumber\\
&& c_k \mathbf{1}=0,\ k\geqslant 2;\quad b_k\mathbf{1}=0, \ k\geqslant -1\}\nonumber\\
\Lambda'&=&\mathbb{C}\{\chi_{-n_1}\dots \chi_{-n_k}\psi_{-m_1}\dots \psi_{-m_\ell} \mathbf{1};\nonumber\\
&& \psi_k \mathbf{1}=0,\ k\geqslant 1;\quad \chi_k\mathbf{1}=0, \ k\geqslant 0\}
\end{eqnarray}
Let's denote $M=\Lambda\otimes \Lambda'$. 
Each of $\Lambda, \Lambda'$ and therefore $M$ has a VOA structure on it, namely, one can define four quantum fields:
\begin{eqnarray}
&&b(z)=\sum_m b_mz^{-m-2}, \qquad c(z)=\sum_n c_nz^{-n+1},\nonumber\\
&&\psi(z)=\sum_m \psi_m z^{-m},\qquad \chi(z)=\sum_m b_mz^{-m-1}
\end{eqnarray}
which according to the commutation relations between modes have the following operator products:
\begin{equation}
b(z)c(w)\sim\frac{1}{z-w},\quad \chi(z)\psi(w)\sim\frac{1}{z-w},
\end{equation}
such that all other operator products do not contain singular terms. 
The Virasoro element is given by the following expression:
\begin{equation}
L^{M}(z)=2:\partial b(z)c(z):+:b(z)\partial c(z):+:\partial \psi(z)\chi(z):,
\end{equation}
such that $b(z)$, $c(z)$ have conformal weights $2$,$-1$, and   
$\psi(z)$, $\chi(z)$ have conformal weights $0$,$1$ correspondingly.
The central charge of the corresponding Virasoro algebra is equal to -28.
One can define the following operator:
\begin{equation}
N_g(z)=:c(z)b(z):+:\psi(z)\chi(z):,
\end{equation}
which is known as ghost number current. 

One can show that the module $W_{\Delta(0),0}^{\varkappa,\eta}$ has the structure of the VOA generate by the quantum fields $\cl(z)=\sum_n\cl_n z^{-n-2}$, and $a(z)=\sum_n a_nz^{-n-1}$, which have the following operator  products (which are equivalent to commutation relations (\ref{w2})):
\begin{eqnarray}
&&\cl(z)\cl(w)\sim \frac{c}{2(z-w)^4}+\frac{2\cl(w)}{(z-w)^2}+\frac{\partial \cl(w)}{z-w},\nonumber\\
&&\cl(z)a(w)\sim \frac{a(w)}{(z-w)^2}+\frac{\partial a(w)}{z-w},\nonumber\\
&& a(z)a(w)\sim \frac{2\eta}{(z-w)^2}.
\end{eqnarray}
Let the space $\mathcal{W}$ be such that $\eta=0, c=28$. Then  the following Proposition is true.
\begin{prop}
The operator of ghost number 1
\begin{eqnarray}
&&Q=\int \frac{\ud z}{2\pi i}J(z), \\
&& J(z)=:c(z)\cl(z):+:c\p cb(z):
\frac{3}{2}\partial^2c(z)+
\psi(z)a(z)+:c\partial\psi\chi(z):\nonumber
\end{eqnarray}
is nilpotent: $Q^2=0$ on $\mathcal{W}\otimes M$.
\end{prop} 
  The space $\mathcal{W}\otimes M$ is known as a semi-infinite cohomology complex, where the differential is $Q$, which is sometimes called the BRST operator. 
The grading in the complex is given by ghost number operator $N_g$. The $k$-th cohomology group
is usually denoted as $H^{\frac{\infty}{2}+k}(W_2,\mathbb{C}{\bf{c}}, \mathcal{W})$.

In this article we will be interested in computing semi-infinite cohomology for 
 $W_2$-modules which have the following form: ${\cal W}=W\otimes 
 \bar {W}$. Here each of $W,\bar{W}$ are representations of $W_2$, with central charges $(c,\eta)$ $(\bar{c},\bar{\eta})$ correspondingly, such that the following relation is satisfied:
 $c+\bar c=28$ and $\eta+\bar \eta=0$. Let us denote the quantum fields for Virasoro  and Heisenberg algebras in $W$ and $\bar {W}$ as $\mathcal{L}(z), a(z)$ and $\bar {\mathcal{L}}(z), {\bar a}(z)$ correspondingly. Then the BRST operator on $\cal W$ has the following form:
 \begin{eqnarray}
&&Q_{\cal W}=\int \frac{\ud z}{2\pi i}J_{\cal W}(z), \\
&& J_{\cal W}(z)=:c(z)(\cl(z)+\bar{\cl}(z)):+:c\p cb(z):+
\frac{3}{2}\partial^2c(z)+\nonumber\\
&&\psi(z)a^{+}(z)+:c\partial\psi\chi(z):\nonumber
\end{eqnarray}
where $a^{+}(z)=a(z)+\bar{a}(z)$. It makes sense to define also the field 
$a^{-}(z)=a(z)-\bar{a}(z)$, which is crucial in the following statement.
\begin{prop}
The operator $Q_{\cal W}$ can be rewritten in the following form:
\begin{eqnarray}
&&Q_{\cal W}=\int \frac{\ud z}{2\pi i}J^{\oplus}_{\cal W}(z), \\
&& J^{\oplus}_{\cal W}(z)=:\tilde c(z)(L(z)+\bar{L}(z)):+:\tilde c\p \tilde c\tilde b(z):+ 
\frac{3}{2}\partial^2\tilde c(z)+
\tilde \psi(z)\tilde a^{+}(z),\nonumber
\end{eqnarray}
Here $L(z)=\cl(z)-\frac{1}{4\eta}:a(z)^2:$, $\bar L(z)=\bar{ \cl}(z)-\frac{1}{4\eta}:\bar a(z)^2:$
and $\tilde A(z)=e^{R}A(z)e^{-R}$ (A is any quantum field),  where $R=\frac{1}{8\eta\pi i }\int dz c(z) a^-(z)\chi(z)$. 
\end{prop}
\noindent {\bf Proof.} 
Let us write explicitly the action of the transformation $\tilde A(z)=e^{R}A(z)e^{-R}$ for every quantum field:
\begin{eqnarray}
&&\tilde a^{+}(z)=a^{+}(z)+\p (c(z)\chi(z)), \quad \tilde b(z)=b(z)-\frac{1}{4\eta}a^-(z)\chi(z),\nonumber\\ 
&&\tilde \psi(z)=\psi(z) +\frac{1}{4\eta}c(z) a^-(z), \quad \tilde a^{-}=a^{-}, \quad 
\tilde c=c, \quad \tilde \chi=\chi.
\end{eqnarray}
To prove this Proposition we just expand the tilded quantum fields in terms of usual ones:
\begin{eqnarray}
&&J^{\oplus}_{\cal W}(z)=: c(z)(L(z)+\bar{L}(z)):+:c\p c\tilde b(z):+ 
\frac{3}{2}\partial^2 c(z)+
\tilde \psi(z)\tilde a^{+}(z)=\nonumber\\
&&: c(z)(L(z)+\bar{L}(z)):+:c(z)\p c(z)(b(z)-\frac{1}{4\eta}a^-(z)\chi(z)):+\frac{3}{2}\partial^2 c(z)+\nonumber\\
&&:(\psi(z)+\frac{1}{4\eta}c(z) a^-(z))(a^{+}(z)+\p (c(z)\chi(z)):=\nonumber\\
&&: c(z)(L(z)+\bar{L}(z)+\frac{1}{4\eta}(a^2(z)-\bar a^2(z))
):+:c\p c\tilde b(z):+ 
\frac{3}{2}\partial^2 c(z)+\nonumber\\
&&\tilde \psi(z)\tilde a^{+}(z)+:\psi(z)\p (c(z)\chi(z)):
\end{eqnarray}

\hfill$\blacksquare$\bigskip

In other words, by means of similarity transformation, one can transform the BRST operator associated with semidirect sum of Virasoro and Heisenberg algebras to the one associated with the direct sum of those algebras. This Proposition has the following Corollary, which will be crucial for computing the semi-infinite cohomology.
\begin{col}
All nontrival cycles of the semi-infinite cohomology belong to the kernels of the following operators: $(i) L_0+\bar L_0+\tilde L^M_0$, $(ii) a_0+\bar{a}_0$, $(iii) T_0$, 
where $T(z)=\frac{1}{4\eta}:\tilde a^{+}(z) \tilde a^{-}(z):+:\p \tilde \psi(z) \tilde\chi(z)$.
\end{col}
\noindent {\bf Proof.} Let us prove $(i)$ first. 
We know that $[Q_{\cal W},\tilde b_0]=L_0+\bar L_0+\tilde L^M_0$, therefore,
\begin{eqnarray}
[Q_{\cal W},\tilde b_0]\Phi=Q_{\cal W}\tilde b_0\Phi=\Delta \Phi.
\end{eqnarray}
Hence, for $\Delta\neq 0$, $\Psi=\Delta^{-1}b_0\Phi$. The proofs of (ii) and (iii) 
are identical to (i), one just needs to use 
the conditions that  $[Q_{\cal W},\tilde \chi_0]=a_0+\bar{a}_0$ and 
$[Q_{\cal W},S_0]=T_0$, where $S_0$ is the zero mode of an operator $S(z)=\frac{1}{4 \eta}a^-(z)\chi(z)$.
\hfill$\blacksquare$\bigskip

\section{ Intertwiners and braided VOAs}

In this section we will study braided VOAs. This will provide a necessary framework 
for the next section, where we will study semi-infinite cohomology and the Lian-Zuckerman multiplication. \\

\noindent {\bf 4.1. Definition of a braided VOA and the simplest example.} First, we give a definition of braided VOAs.\\ 

\noindent{\bf Definition.}\cite{fz} {\it Let $\mathbb{V}=\oplus_{\lambda\in I}\mathbb{V}_{\lambda}$ be a direct sum of graded complex 
vector spaces, called sectors: $\mathbb{V}_{\lambda}=\oplus_{n\in \mathbb{Z}_+}\mathbb{V}_{\lambda}[n]$, indexed by some set $I$. Let $\Delta_{\lambda}$, $\lambda\in I$ be complex numbers, which we will call conformal weights of the corresponding sectors.
We say that $\mathbb{V}$ is a braided vertex algebra, if there are distinguished elements $0\in I$ such  that $\Delta_0=0$, $\mathbf{1}\in \mathbb{V}_0[0]$, linear maps $D:\mathbb{V}\to\mathbb{V}$, $\mathcal{R}:\mathbb{V}\otimes\mathbb{V}\to \mathbb{V}\otimes\mathbb{V}$ and the linear correspondence 
\begin{eqnarray}
\mathbb{Y}(\cdot,z)\cdot:\mathbb{V}\otimes\mathbb{V}\to \mathbb{V}\{z\},\quad \mathbb{Y}=\sum_{\lambda,\lambda_1,\lambda_2}\mathbb{Y}^{\lambda_1\lambda_2}_{\lambda}(z),
\end{eqnarray}
where 
\begin{eqnarray}
\mathbb{Y}^{\lambda_1\lambda_2}_{\lambda}(z)\in Hom(\mathbb{V}_{\lambda_1}\otimes\mathbb{V}_{\lambda_2},\mathbb{V}_{\lambda})\otimes z^{\Delta_{\lambda}-\Delta_{\lambda_1}-\Delta_{\lambda_2}}
\mathbb{C}[[z,z^{-1}]],
\end{eqnarray}
such that the following properties are satisfied:\\
i)Vacuum property: $\mathbb{Y}(\mathbf{1},z)v=v$, $\mathbb{Y}(v,z)\mathbf{1}|_{z=0}=v$.\\
ii) Complex analyticity: for any $v_i\in \mathbb{V}_{\lambda_i}$, $(i=1,2,3,4)$ the matrix elements 
$\langle v_4^*, \mathbb{Y}(v_3,z_2)\mathbb{Y}(v_2,z_1)v_1 \rangle$ regarded as formal Laurent series in $z_1,z_2$, converge in the domain $|z_2|> |z_1|$ to a complex analytic function $r(z_1,z_2)\in z_1^{h_1}z_2^{h_2}
(z_1-z_2)^{h_3} \mathbb{C}[z_1^{\pm 1},z_2^{\pm 1}, (z_1-z_2)^{-1}]$, where 
$h_1, h_2,h_3\in \mathbb{C}$.\\
iii) Derivation property: $\mathbb{Y}(Dv,z)\mathbf{1}=\frac{d}{dz}\mathbb{Y}(v,z).$\\
iv) Braided commutativity (understood in the weak sense \footnote{ By the $weak$ sense we mean that the relation holds for the matrix elements of the corresponding operator products. }):
\begin{eqnarray}
\mathscr{A}_{z,w}(\mathbb{Y}(v,z)\mathbb{Y}(u,w))=\sum_i\mathbb{Y}(u_i,w)\mathbb{Y}(v_i,z),
\end{eqnarray}
where $\mathcal{R}(u\otimes v)=\sum_i u_i\otimes v_i$ and 
 $\mathscr{A}_{z,\omega}$ denote the monodromy around the path
\begin{eqnarray}
w(t)&=&\frac{1}{2}\big((z+w)+(w-z)e^{\pi it}\big),\nonumber\\
     z(t)&=&\frac{1}{2}\big((z+w)+(z-w)e^{\pi it}\big),
\end{eqnarray}
as shown on the picture:
\vspace{3mm}
\begin{eqnarray}
\xymatrix{
 z\ \bullet\ \ar@/_1pc/[r] &
\bullet \   \ar@/_1pc/[l] w
}
\end{eqnarray}

\vspace{3mm}

\noindent v) There exists an element $\omega\in \mathbb{V}_0$, such that
\begin{equation}
Y(\omega,z)=\sum_{n\in\mathbb{Z}}L_nz^{-n-2}
\end{equation}
and $L_n$ satisfy the relations of Virasoro algebra with $L_{-1}=D$.\\
vi) Associativity (understood in the weak sense):
\begin{eqnarray}
\mathbb{Y}(\mathbb{Y}(u,z-w)v, w)=\mathbb{Y}(u,z)\mathbb{Y}(v,w).
\end{eqnarray}}
The associativity condition puts a restriction on the operator $\mathcal{R}$, namely one can show that 
matrix $R$ satisfies the Yang-Baxter equation (see \cite{fz}, \cite{styrkas}).

\noindent The simplest example of braided VOA (which does not reduce to usual VOA) is a natural extension of the VOA $F_{\eta}(a)$ generated by Heisenberg algebra, considered in the previous section. Let us consider the following space:
\begin{equation}
\hat{F}_\eta=\oplus_{\lambda\in\mathbb{Z\eta}}F_{\lambda,\eta}.
\end{equation}
Below we will show that $\hat{F}_\eta$ carries a structure of braided VOA.
The operators
\begin{equation}\label{xz}
\mathbb{X}(\lambda,z)=\mathbf{1}_\lambda z^{\frac{\lambda a_0}{2\eta}}
e^{\big(\frac{\lambda}{2\varkappa}\sum_{n>0}\frac{a_{-n}}{n}z^n\big)}
e^{-\big(\frac{\lambda}{2\varkappa}\sum_{n>0}\frac{a_n}{n}z^{-n}\big)},
\end{equation}
where $\lambda\in\mathbb{Z}$, generate vacuum vectors in Fock modules. 
It is clear that $\mathbb{X}(\lambda,z)\mathbf{1}_0|_{z=0}=\mathbf{1}_\lambda$. 
Denoting
\begin{equation}
\mathbb{X}_{n_1,\dots,n_k}(\lambda,z)\equiv :a^{(n_1)}(z)\dots a^{(n_k)}(z)\mathbb{X}(\lambda,z):,
\end{equation}
where $a^{(n)}(z)=\frac{1}{n!}\big(\frac{\ud}{\ud z}\big)^na(z)$, one
can see that
\begin{equation}
\mathbb{X}_{n_1,\dots,n_k}(\lambda,z)\mathbf{1}_0|_{z=0}=a_{-n_1},\dots,a_{-n_k}\mathbf{1}_\lambda.
\end{equation}
In such a way, we build the correspondence
\begin{equation}
i:v\rightarrow Y(v,z)=\sum_{n\in \mathbb{Z}}v_{(n)}z^{-n-1},
\end{equation}
such that $v\in\hat{F}_\eta$ and $v_{(n)}\in \mathrm{End}(\hat{F}_\eta)$. 
To see the connection with definition of braided VOA, we note that in this case the set $I=\mathbb{Z}$, the sectors are Fock spaces $F_{\lambda,\eta}$, 
$\Delta(\lambda)$ is the conformal weight of the vacuum.

\noindent Let $|z|>|w|$, then
\begin{equation}
\mathbb{X}(\lambda,z)\mathbb{X}(\mu,w)=(z-w)^{\frac{\lambda\mu\eta}{2}}
(\mathbb{X}(\lambda+\mu,w)+\dots),
\end{equation}
where dots stand for the terms regular in $(z-w)$, hence
\begin{eqnarray}\label{acon}
\mathscr{A}_{z,w}\big(\mathbb{X}(\lambda,z)\mathbb{X}(\mu,w)\big)=
p^{\frac{\lambda\mu}{2}}\mathbb{X}(\mu,w)\mathbb{X}(\lambda,z),
\end{eqnarray}
where $p=e^{\pi i\eta}$. Here, we underline that 
the expression above should be understood in a weak sense, i.e. the analytical continuation is performed for the matrix elements of the corresponding operator products. Moreover, the matrix elements of operator product expansion 
$\mathbb{X}(\lambda,z)\mathbb{X}(\mu,w)$ exist in the domain $|z|>|w|$ and the analytical continuation relates it to the matrix elements of operator product expansion $\mathbb{X}(\mu,w)\mathbb{X}(\lambda,z)$, which converge in the domain $|w|>|z|$.
The relation (\ref{acon}) is a simplified case of associativity condition, 
since in our case $\mathbb{Y}^{\lambda}_{\lambda_1\lambda_2}$ is nonzero only for $\lambda=\lambda_1+\lambda_2$ and the $\mathcal{R}$-operator is therefore reduced to the multiplication on some power of $p$.

In the remainig part of this section we will discuss examples which are more involved. In order to construct them we need to use additional constructions.\\

\noindent{\bf 4.3. Intertwiners for Virasoro and $W_2$ algebras.} 
In \cite{fz} we have constructed an intertwiner for Virasoro algebra. This is a 
map 
\begin{equation}\label{intvir}
\Phi_{\lambda\mu}^\nu(z):V_{\Delta(\lambda),\varkappa}\otimes V_{\Delta(\mu),\varkappa}\rightarrow 
V_{\Delta(\nu),\varkappa}[[z,z^{-1}]]z^{\Delta(\nu)-\Delta(\mu)-\Delta(\lambda)},
\end{equation}
which has the following property \cite{ms}:
\begin{equation}\label{propint}
L_n\cdot\Phi_{\lambda\mu}^\nu(z)=\Phi_{\lambda\mu}^\nu(z)\Delta_{z,0}(L_n),
\end{equation}
where
\begin{equation}
\Delta_{z,0}(L_n)=\oint_z\frac{\ud\xi}{2\pi i}\xi^{n+1}\Big(\sum_m(\xi-z)^{-m-2}L_m\Big)\otimes 1+1\otimes L_n.
\end{equation}
Moreover, we have found that they satisfy the relation involving the braiding matrix $B^V$ of $U_q(sl(2))$.
\begin{prop}{\rm \cite{fz}}
Let $z_1,z_2\in \mathbb{C}$ such that $0<|z_1|<|z_2|$, $\lambda_i\ge 0$ $(i=0,1,2,3)$. Then the following relation holds:
\begin{eqnarray}\label{fockint}
&&\mathscr{A}_{z_1,z_2}\big(\Phi_{\lambda_3\rho}^{\lambda_0}(z_2)\Phi_{\lambda_2\lambda_1}^\rho(z_1)\big)(P\otimes 1)=\nonumber\\
&&\sum_{\xi}B^V_{\rho\xi}
\left[\begin{array}{cc}
\lambda_0 &\lambda_1 \\
\lambda_2 &\lambda_3 
\end{array}
\right]
\Phi_{\lambda_2\xi}^{\lambda_0}(z_1)\Phi_{\lambda_3\lambda_1}^\xi(z_2), 
\end{eqnarray} 
where $P$ is an interchange operator, namely $P(v_1\otimes v_2)=v_2\otimes v_1$ and $q=e^{\frac{\pi i}{\varkappa}}$.
\end{prop}

Now we show that there exists an intertwining operator between the corresponding irreducible highest weight representations of $W_2$ algebra.
\begin{prop}
i)There exists a map
\begin{eqnarray}\label{intw2}
&&\Phi_{\lambda,l;\mu,m}^{\nu,n}(z):W^{\eta,\kappa}_{\Delta(\lambda),l}\otimes W^{\eta,\kappa}_{\Delta(\mu),m}\rightarrow \nonumber\\
&&W^{\eta,\kappa}_{\Delta(\nu),n}[[z,z^{-1}]]z^{\eta(n^2-m^2-l^2)}z^{\Delta(\nu)-\Delta(\mu)-\Delta(\lambda)},
\end{eqnarray}
such that 
\begin{eqnarray}
&&\cl_n\cdot\Phi_{\lambda,l;\mu,m}^{\nu,n}(z)=\Phi_{\lambda,l;\mu,m}^{\nu,n}(z)\Delta_{z,0}(\cl_n),\nonumber\\
&&a_k\cdot\Phi_{\lambda,l;\mu,m}^{\nu,n}(z)=\Phi_{\lambda,l;\mu,m}^{\nu,n}(z)\Delta_{z,0}(a_k),
\end{eqnarray}
where
\begin{eqnarray}
&&\Delta_{z,0}(\cl_n)=\oint_z\frac{\ud\xi}{2\pi i}\xi^{n+1}\Big(\sum_m(\xi-z)^{-m-2}\cl_m\Big)\otimes 1+1\otimes \cl_n ,\nonumber\\
&&\Delta_{z,0}(a_k)=\oint_z\frac{\ud\xi}{2\pi i}\xi^{k}\Big(\sum_m(\xi-z)^{-m-1}a_m\Big)\otimes 1+1\otimes a_k.
\end{eqnarray}
ii) The intertwining operators  $\Phi_{\lambda,l;\mu,m}^{\nu,n}(z)$ satisfy the following quadratic relation:
\begin{eqnarray}\label{w2intrel}
&&\mathscr{A}_{z_1,z_2}\big(\Phi_{\lambda_3,l_3;\rho,r}^{\lambda_0,l_0}(z_2)\Phi_{\lambda_2,l_2;\lambda_1,l_1}^{\rho,r}(z_1)\big)(P\otimes 1)=\nonumber\\
&&\Big(\frac{p}{q}\Big)^{l_2l_3}\sum_{\xi,k}B_{\rho,r;\xi,k}
\left[\begin{array}{cc}
\lambda_0,l_0 &\lambda_1,l_1 \\
\lambda_2,l_2 &\lambda_3,l_3 
\end{array}
\right]
\Phi_{\lambda_2,l_2;\xi,k}^{\lambda_0,l_0}(z_1)\Phi_{\lambda_3,l_3;\lambda_1,l_1}^{\xi,k}(z_2), 
\end{eqnarray} 
where $p=e^{\pi i\eta}$, $q=e^{\frac{\pi i}{\varkappa}}$.
\end{prop}
\noindent {\bf Proof.} i) First, we construct an intertwining operator
\begin{eqnarray}
\Psi^n_{lm}:F_{l,\eta}\otimes F_{m,\eta}\to F_{n,\eta}[[z,z^{-1}]]z^{\eta(n^2-m^2-l^2)},
\end{eqnarray}
which obeys the rule 
\begin{eqnarray}\label{aint}
a_k\cdot\Psi_{lm}^{n}(z)=\Phi_{lm}^{n}(z)\Delta_{z,0}(a_k).
\end{eqnarray}
However, we already constructed one, when studied the first 
example of braided algebra in subsection 4.1. The operator
\begin{eqnarray}
\delta_{l+m,n}Y(\cdot, z)\cdot:F_{l,\eta}\otimes F_{m,\eta}\to F_{n,\eta}[[z,z^{-1}]]z^{\eta(n^2-m^2-l^2)}
\end{eqnarray}
satisfies the property (\ref{aint}). Now we introduce an operator 
\begin{eqnarray}
\Phi_{\lambda,l;\mu,m}^{\nu,n}(z)=\Phi_{\lambda,\mu}^{\nu}(z)\otimes \Psi_{lm}^{n}(z),
\end{eqnarray}
where $\Phi_{\lambda,\mu}^{\nu}(z)$ is the intertwining operator for irreducible highest weight modules of Virasoro algebra, generated by $L_n\equiv \cl_n-\frac{1}{4\eta}\sum_m:a_{n-m}a_m:$ operators, and therefore is the intertwining operator for the irreducible highest weight representations of $W_2$ algebra. 

The second statement of the Proposition follows from the commutativity condition of $Y$ and similar statement for intertwining operators for Virasoro algebra.
\hfill$\blacksquare$\bigskip

\noindent{\bf 4.4. From $W_2$ and $U_q(gl(2))$ to braided VOA.} 
Now we have all necessary tools to build a more sophisticated example of a 
braided vertex algebra than the one in section 4.1. In \cite{fz} we have constructed the braided VOA 
on the space 
\begin{equation}
\mathbb{F}_\varkappa =\bigoplus_{\lambda\in\mathbb{Z}_+} (V_{\Delta(\lambda),\varkappa}\otimes V_\lambda), 
\end{equation}
where 
$V_{\lambda}$ is an irreducible representation of $U_q(sl(2))$ with highest weight 
$\lambda$ and $V_{\Delta(\lambda)}$ are highest weight Virasoro modules discussed in Subsection 3.1., such that the $\mathcal{R}$-operator from braided commutativity relation of $\mathbb{F}_{\varkappa}$ is related to the universal $R$-matrix for $U_q(sl(2))$.

Here we modify the construction of \cite{fz} in order to build braided VOA compatible with the $W_2$ and $U_q(gl(2))$ structures, namely, it will be defined on the space 
\begin{equation}
\mathbb{G}_{\varkappa,\eta}= \bigoplus_{\lambda\in\mathbb{Z}_+, k\in\mathbb{Z}}(W^{\varkappa,\eta}_{\Delta(\lambda), k}\otimes V_{\lambda,k}).
\end{equation}
We define a map
\begin{equation}
Y(\cdot, z)\cdot: \mathbb{G}_{\varkappa,\eta}\otimes \mathbb{G}_{\varkappa,\eta}\to \mathbb{G}_{\varkappa,\eta}\{z\}
\end{equation}
in such a way that 
\begin{equation}\label{bvoa}
Y:v\otimes a\rightarrow Y(v\otimes a,z)=
\sum_{\nu,\mu} \Phi^{\nu,n}_{\lambda,l;\mu,m}(z)(v\otimes \cdot)\otimes \phi^{\nu,n}_{\mu,m;\lambda,l}(\cdot \otimes a),
\end{equation}
where $v\in W^{\varkappa,\eta}_{\Delta(\lambda),l}$ and $a\in V_{\lambda,l}$.

Let us denote $\mathbb{G}_{\varkappa,\eta}(l,\lambda)\equiv W^{\varkappa,\eta}_{\Delta(\lambda), l}\otimes V_{\lambda,l}$.  One can see that the sectors for the map $Y$ are given by the spaces $\mathbb{G}_{\varkappa,\eta}(l,\lambda)$ and the conformal weight of the sector is $\Delta(\lambda)+\eta k^2$. Let us prove that $Y$ satisfies the braided commutativity relation and compute an explicit expression for the $\mathcal{R}$-operator.
Let $v_i\otimes a_i\in  W^{\varkappa,\eta}_{\Delta(\lambda_i), l_i}\otimes V_{\lambda_i,l_i}$ $(i=1,2,3)$, then
\begin{eqnarray}
&&\mathscr{A}_{z_2,z_1}\big(Y(v_1\otimes a_1,z_2)Y(v_2\otimes a_2,z_1)\big)(v_3\otimes a_3)=\nonumber\\
&&\mathscr{A}_{z_2,z_1}\Big(\sum_{\lambda_1,\lambda_2,\nu,\rho}
\Phi_{\lambda_1,l_1\rho,r}^{\nu,n}(z_2)\Phi_{\lambda_2,l_2;\lambda_3,l_3}^{\rho,r}(z_1)\otimes
\phi_{\rho,r;\lambda_1,l_1}^{\nu,n}\phi_{\lambda_3,l_3;\lambda_2,l_2}^{\rho,r}\Big)\cdot\nonumber\\
&&(v_1\otimes v_2\otimes v_3)\otimes(a_3\otimes a_2\otimes a_1)=\nonumber\\
&&\Big(\sum_{\lambda_1,l_1,\lambda_2,l_2,\nu,n,\rho,r,\xi,k}\Phi_{\lambda_2,l_2,\xi,k}^{\nu,n}(z_1)\Phi_{\lambda_1,l_1;\lambda_3,l_3}^{\xi,k}(z_2)\nonumber\\
&&
\Big(\frac{p}{q}\Big)^{l_1l_2}
B_{\rho,r;\xi,k}
\left[\begin{array}{cc}
\nu,n &\lambda_3,l_3 \\
\lambda_2,l_2 &\lambda_1,l_1 
\end{array}
\right] \otimes 
\phi_{\rho,r;\lambda_1,l_1}^{\nu,n}\phi_{\lambda_3,l_3;\lambda_2,l_2}^{\rho,r}\Big)\cdot\nonumber\\
&&(v_2\otimes v_1\otimes v_3)(a_3\otimes a_2\otimes a_1)=\nonumber\\
&&\Big(\sum_{\lambda_1,l_1,\lambda_2,l_2,\nu,n,\xi,k}
\Phi_{\lambda_2,l_2;\xi,k}^{\nu,n}(z_1)\Phi_{\lambda_2,l_2;\lambda_3,l_3}^{\xi,k}(z_2)\otimes
\phi_{\xi,k;\lambda_2,l_2}^{\nu,n}\phi_{\lambda_3,l_3;\lambda_1,l_1}^{\xi,k}\Big)\cdot\nonumber\\
&&(v_2\otimes v_1\otimes v_3)\otimes(a_3\otimes \sum_ir_i^{(2)}a_1\otimes r_i^{(1)}a_2)=\nonumber\\
&&\sum_i Y(v_2\otimes r_i^{(1)}a_2,z_1)Y(v_1\otimes r_i^{(2)}a_1,z_1)(v_3\otimes a_3),
\end{eqnarray}
where
\begin{equation} 
\mathcal{R}\equiv r^{(1)}_i\otimes r_i^{(2)}=\Big(\frac{p}{q}\Big)^{I\otimes I}R, 
\end{equation}
and $R$ is a universal R-matrix for $U_q(gl(2))$. Therefore $Y$ satisfies the braided commutativity relation. One can prove that $Y$ satisfies all other necessary properties along the lines of \cite{fz} and we arrive to the following Proposition.
\begin{prop}
The map $Y$ defines a structure of braided VOA on 
$\mathbb{G}_{\varkappa,\eta}$.
\end{prop}
The braided VOA $\mathbb{G}_{\varkappa,\eta}$ possesses the following remarkable subalgebras.
\begin{col}
$\mathbb{G}_{\varkappa,\eta}$ possesses braided vertex subalgebras 
$\hat{\mathbb{G}}_{\varkappa,\eta}$, 
 $\hat{\mathbb{G}^{+}}_{\varkappa,\eta}$ with the spaces
\begin{eqnarray}
&&\hat{\mathbb{G}}_{\varkappa,\eta}=\bigoplus_{\lambda\in\mathbb{Z}_+, r\in\mathbb{Z}}(W^{\varkappa,\eta}_{\Delta(\lambda), \lambda+2r}\otimes V_{\lambda,\lambda+2r}),\nonumber\\
&&\hat{\mathbb{G}^{+}}_{\varkappa,\eta}=\bigoplus_{\lambda,r\in\mathbb{Z}_+}(W^{\varkappa,\eta}_{\Delta(\lambda), \lambda+2r}\otimes V_{\lambda,\lambda+2r}).
\end{eqnarray}
Here $\hat{\mathbb{G}}^+_{\varkappa,\eta}$ is the minimal braided subalgebra of $\mathbb{G}_{\varkappa,\eta}$, 
containing subspace $\mathbb{G}_{\varkappa,\eta}(1,1)$, and  $\hat{\mathbb{G}}_{\varkappa,\eta}$ can be obtained from it by extending via subspace $\mathbb{G}_{\varkappa,\eta}(-2,-2)$.
\end{col}
Finally, we mention the following important property. There is a natural action of 
$U_q(gl(2))$ on ${\mathbb{G}}_{\varkappa,\eta},\ \hat{\mathbb{G}}_{\varkappa,\eta}, \hat{\mathbb{G}}^+_{\varkappa,\eta}$. The next Proposition shows that it is compatible with the braided VOA structure.
\begin{prop}
There is a natural $U_q(gl(2))$ action on the vertex algebras 
${\mathbb{G}}_{\varkappa,\eta},\ \hat{\mathbb{G}}_{\varkappa,\eta}, \hat{\mathbb{G}}^+_{\varkappa,\eta}$, such that $gY=Y\Delta(g),$
where $g\in U_q(gl(2))$. 
\end{prop}

\noindent {\bf Proof.} This is a direct consequence of the definition of $Y$ and the properties of the intertwining operator on $U_q(gl(2))$. 
\hfill$\blacksquare$\bigskip

\section{ $GL_q(2)$ as semi-infinite cohomology and Laplace operator.}
{\bf 5.1. Computation of the semi-infinite cohomology.} We are interested in computation of the semi-infinite cohomology of the tensor product of the following modules of $W_2$:
\begin{equation}
W^{\varkappa,\eta}_{\Delta(\lambda),l}\otimes W^{-\varkappa,-\eta}_{\Delta(\mu),m}.
\end{equation}
For $H^{\frac{\infty}{2}+0}(W_2,\mathbb{C}{\bf c},\cdot)$ the answer is given in the following proposition.

\begin{prop}
The 0th semi-infinite cohomology group for the tensor product of two irreducible highest weight $W_2$ modules is given by
\begin{equation}
H^{\frac{\infty}{2}+0}
(W_2,\mathbb{C}{\bf c},W^{\varkappa,\eta}_{\Delta(\lambda),l}\otimes W^{-\varkappa,-\eta}_{\Delta(\mu),m})=\mathbb{C}\delta_{\lambda,\mu}\delta_{l,m}.
\end{equation}
\end{prop}
\noindent {\bf Proof.} First, to simplify the semi-infinite cohomology operator we use Proposition 3.5. 
As well as there we have tensor product of $W_2$-modules with
complimentary central charges. Moreover, the associated semi-infinite complexes formed out of "tilded" variables and original ones are isomorphic to each other. Therefore, we can compute the semi-infinite cohomology using the operator $Q^{\oplus}$. Part (ii) of Corollary 3.2. says that $H^{\frac{\infty}{2}+0}$ is nontrivial iff $l=m$. Part (i) says that $\lambda$ should be equal to $\mu$, since $\varkappa$ is generic. Part (iii) leads to the fact that only the highest weight vectors of the tensor product of Fock modules $F_{\varkappa,\eta}\otimes F_{-\varkappa,-\eta}\otimes M_{\psi,\chi}$ (where $M_{\psi,\chi}$ is the part of $M$ generated by $\psi_m,\chi_n$ modes) contribute to nontrivial cohomology classes. Hence, the problem is reduced  to the computation of the 0th semi-infinite cohomology group of Virasoro modules, i.e. 
$H^{\frac{\infty}{2}+0}(Vir,\mathbb{C}{\bf c},V_{\Delta(\lambda),\varkappa}\otimes V_{\Delta(\lambda), -\varkappa})$, which is known to be equal to $\mathbb{C}$ (see \cite{lz1}, \cite{lz2}, \cite{fs}). Hence, the proposition is proven.
\hfill$\blacksquare$\bigskip

\noindent Let us introduce the following spaces with the structure of braided VOA:
\begin{eqnarray}
&&\mathbb{G}=\mathbb{G}_{\varkappa,\eta}\otimes \mathbb{G}_{-\varkappa,-\eta}, \quad \hat{\mathbb{G}}=\hat{\mathbb{G}}_{\varkappa,\eta}\otimes \hat{\mathbb{G}}_{-\varkappa,-\eta}\nonumber\\
&&\hat{\mathbb{G}}^+=\hat{\mathbb{G}}^+_{\varkappa,\eta}\otimes \hat{\mathbb{G}}^+_{-\varkappa,-\eta}
\end{eqnarray}
We see that $\hat{\mathbb{G}}$ is a braided VOA subalgebra in $\mathbb{G}$.
An immediate consequence of Proposition 5.1. is the following Theorem.
\begin{theorem}
The 0th cohomology groups of the spaces of braided vertex algebras 
$\mathbb{G}$, $\hat{\mathbb{G}}$,  $\hat{\mathbb{G}}^+$ are given by 
\begin{eqnarray}
&& H^{\frac{\infty}{2}+0}(Vir,\mathbb{C}{\bf c},\mathbb{G})=
\oplus_{\lambda,k\in\mathbb{Z}}V_{\lambda,k}\otimes \bar{V}_{\lambda,k},\nonumber\\
&& H^{\frac{\infty}{2}+0}(Vir,\mathbb{C}{\bf c},\hat{\mathbb{G}})=
\oplus_{\lambda\in\mathbb{Z}_+,r\in\mathbb{Z}}V_{\lambda,\lambda+2r}\otimes \bar{V}_{\lambda,\lambda+2r},\nonumber\\
&& H^{\frac{\infty}{2}+0}(Vir,\mathbb{C}{\bf c},\hat{\mathbb{G}}^+)=
\oplus_{\lambda,r\in\mathbb{Z}_+}V_{\lambda,\lambda+2r}\otimes \bar{V}_{\lambda,\lambda+2r}
\end{eqnarray}
where $\bar{V}_{\lambda,k}$ stands for irreducible representation of $U_q(gl(2))$ after involution $q \to q^{-1}$.
\end{theorem}

Now we will define a product structure on these cohomology spaces using the natural product structure on VOA.\\

\noindent{\bf 5.2. Ring structure on the semi-infinite cohomology spaces.} The Lian-Zuckerman associative product structure is defined on the representatives $U$, $V$ of the cohomology classes of VOA as follows \cite{lz2}:
\begin{equation}\label{fzprod}
\mu(U,V)=\mathrm{Res}_z \Big(\frac{U(z)V}{z}\Big).
\end{equation}
In \cite{fz} we have shown that it is also gives an associative algebra structure 
on the space $\mathbb{F}=\mathbb{F}_{\varkappa}\otimes\mathbb{F}_{\varkappa}$. The same reasoning we used in that case applies to $\mathbb{G}$. Hence, we obtain the following Proposition.
\begin{prop}
The 0th semi-infinite cohomology space of the braided VOAs $\mathbb{G}$, $\hat {\mathbb{G}}$,  $\hat {\mathbb{G}^+}$ 
possess an associative product given by the formula (\ref{fzprod}) on the representatives of cohomology classes, such that it has the following subalgebras:
\begin{eqnarray}
&&(H^{\frac{\infty}{2}+0}(W_2,\mathbb{C}\mathbf{c},\hat{\mathbb{G}}^+),\mu)\subset\nonumber\\
&&(H^{\frac{\infty}{2}+0}(W_2,\mathbb{C}\mathbf{c},\hat{\mathbb{G}}),\mu)
\subset (H^{\frac{\infty}{2}+0}(W_2,\mathbb{C}\mathbf{c},\mathbb{G}),\mu).
\end{eqnarray} 
 \end{prop}

\noindent This proposition has an immediate corollary, which follows from the braided commutativity property.

\begin{col}
The operation $\mu$ being considered on $H^{\frac{\infty}{2}+0}(W_2,\mathbb{C}\mathbf{c},\mathbb{G})$ is associative and satisfies the following commutativity relation:
\begin{equation}
\mu(U,V)=\mu(\hat{r}_i^{(1)}V,\hat{r}_i^{(2)}U).
\end{equation}
Here $\hat{\mathcal{R}}=\sum_i \hat{r}_i^{(1)}\otimes\hat{r}_i^{(2)}=\mathcal{R}\bar{\mathcal{R}}$, where 
$\mathcal{R}$, $\bar{\mathcal{R}}$ are the braiding operators on $\mathbb{G}_{\varkappa}, \mathbb{G}_{-\varkappa}$ correspondingly.  
\end{col}
Now we want to find the generating set for both $H^{\frac{\infty}{2}+0}(W_2,\mathbb{C}\mathbf{c},\hat{\mathbb{G}})$ and $H^{\frac{\infty}{2}+0}(W_2,\mathbb{C}\mathbf{c},\mathbb{G})$ under the multiplication $\mu$. 

We recall that we denoted $\mathbb{G}_{\varkappa,\eta}(l,\lambda)\equiv W^{\varkappa,\eta}_{\Delta(\lambda), l}\otimes V_{\lambda,l}$. Let $\mathbb{G}(\mu,k)=\mathbb{G}_{\varkappa,\eta}(\mu,k)\otimes\mathbb{G}_{-\varkappa,-\eta}(\mu,k)$. Due to the structure of the braided VOA on $\mathbb{G}$ and 
$\hat{\mathbb{G}}$ we have the following Proposition.
\begin{prop}
i) Algebra $(H^{\frac{\infty}{2}+0}(W_2,\mathbb{C}\mathbf{c},\hat{\mathbb{G}}^+),\mu)$ is generated by \\
$H^{\frac{\infty}{2}+0}(W_2,\mathbb{C}\mathbf{c},{\mathbb{G}}(1,1))$.\\
\noindent ii) Algebra $(H^{\frac{\infty}{2}+0}(W_2,\mathbb{C}\mathbf{c},\hat{\mathbb{G}}),\mu)$ is generated by \\
$H^{\frac{\infty}{2}+0}(W_2,\mathbb{C}\mathbf{c},{\mathbb{G}}(1,1))$ and $H^{\frac{\infty}{2}+0}(W_2,\mathbb{C}\mathbf{c},{\mathbb{G}}(0,-2))$.\\
\noindent iii)Algebra $(H^{\frac{\infty}{2}+0}(W_2,\mathbb{C}\mathbf{c},{\mathbb{G}}),\mu)$ is generated by $H^{\frac{\infty}{2}+0}(W_2,\mathbb{C}\mathbf{c},{\mathbb{G}}(1,1))$ and $H^{\frac{\infty}{2}+0}(W_2,\mathbb{C}\mathbf{c},{\mathbb{G}}(0,-1))$.
\end{prop}
\noindent We note that $H^{\frac{\infty}{2}+0}(W_2,\mathbb{C}\mathbf{c},\mathbb{G}(1,1))$ $\cong
V_{1,1}\otimes \bar{V}_{1,1}$ and $H^{\frac{\infty}{2}+0}(W_2,\mathbb{C}\mathbf{c},{\mathbb{G}}(0,-2))$ $\cong\mathbb{C}$, $H^{\frac{\infty}{2}+0}(W_2,\mathbb{C}\mathbf{c},\hat{\mathbb{G}}(0,-1))\cong\mathbb{C}$.

Since $\mathcal{R}=p^{I\otimes I}\mathcal{R}'$, where $\mathcal{R}'$ is the universal R-matrix for $U_q(sl(2))$ (i.e. the braiding operator from braided commutativity relation for $\mathbb{F}_{\varkappa}$), we have 
$\hat{\mathcal{R}}=R'\bar{R'}$. 
Now we calculate the commutation relations between the elements from the generating set. 
There are only two vectors in each of 
$V_{1,1},\bar V_{1,1}$, i.e. the highest weight and the lowest weight vectors.
We denote them as $a_+, a_-$ and $\bar{a}_+,\bar{a}_-$ correspondingly. Let us make the following notation:
\begin{eqnarray}\label{abcd}
v\otimes a_-\otimes\bar{a}_+=A,\quad
v\otimes a_+\otimes\bar{a}_-=D,\nonumber\\
v\otimes a_+\otimes\bar{a}_+=B,\quad
v\otimes a_-\otimes\bar{a}_-=C,
\end{eqnarray}
where 
$v\in H^{\frac{\infty}{2}+0}(W_2,\mathbb{C}\mathbf{c},W^{\varkappa,\eta}_{1,1}\otimes W^{-\varkappa,-\eta}_{1,1})$. 
Then the following Proposition holds, which is a consequence of Proposition 6.9 of \cite{fz}.
\begin{theorem}
i)The generators $A,B,C,D$ satisfy the following relations:
\begin{eqnarray}
&&AB=BA{q^{-1}}, \quad CB=BC,\quad DB=BDq,\quad CA=ACq,\nonumber\\
&&AD-DA=(q^{-1}-q)BC,\quad CD=DCq^{-1},
\end{eqnarray}
and $det_q\equiv AD-q^{-1}BC$ is a nonzero element in $H^{\frac{\infty}{2}+0}(W_2,\mathbb{C}\mathbf{c},\hat{\mathbb{G}}(0,2))$.\\
\noindent ii)  The element $det_q$ has inverse, which belongs to $H^{\frac{\infty}{2}+0}(W_2,\mathbb{C}\mathbf{c},\hat{\mathbb{G}}(0,-2))$, i.e.
$(H^{\frac{\infty}{2}+0}(W_2,\mathbb{C}\mathbf{c},\hat{\mathbb{G}}),\mu)\cong GL_q(2)$.
\end{theorem}
\noindent {\bf Remark. } One can see that the algebra  
$(H^{\frac{\infty}{2}+0}(W_2,\mathbb{C}\mathbf{c},\hat{\mathbb{G}}^+),\mu)$ is a subalgebra of $GL_q(2)$ generated by the elements $A,B,C,D$ only, 
and $(H^{\frac{\infty}{2}+0}(W_2,\mathbb{C}\mathbf{c},{\mathbb{G}}),\mu)$
corresponds to $GL_q(2)$ extended by element $t\in H^{\frac{\infty}{2}+0}(W_2,\mathbb{C}\mathbf{c},\hat{\mathbb{G}}(0,1))$, such that $t^2=det_q$. \\

\noindent {\bf 5.3. Quantum Laplacian.}
Let us denote $GL_q(2)^+\equiv (H^{\frac{\infty}{2}+0}(W_2,\mathbb{C}\mathbf{c},\hat{\mathbb{G}}^+),\mu)$. In \cite{fj} it was defined an extension of this algebra by the invertible element $\delta$, obeying the following commutation relations with the generating elements:
\begin{eqnarray}
\delta A=A\delta, \quad \delta B=qB\delta, \nonumber\\
\delta C=q^{-1}C\delta, \quad \delta D=D\delta.
\end{eqnarray}
In such a way one can define elements $x_{11}, x_{12}, x_{21}, x_{22}$:
\begin{eqnarray}
x_{11} = \delta A, \quad x_{12} = q^{-1/2}\delta B ,\nonumber\\
x_{21} = q^{1/2}\delta C ,\quad  x_{22} = \delta D
\end{eqnarray}
satisfying commutation relations:
\begin{eqnarray}
&&x_{11}x_{12} = x_{12}x_{11} ,\quad x_{21}x_{22} = x_{22}x_{21} , \nonumber\\
&&[x_{11} , x_{22} ] + [x_{21} , x_{12} ] = 0 , \nonumber\\
&&x_{11}x_{21} = q^{-2}x_{21}x_{11} , \quad x_{12}x_{22} = q^{-2}x_{22}x_{12} , \nonumber\\
&&x_{21}x_{12} = q^2x_{12}x_{21}.
\end{eqnarray}
The space generated by $x_{11},x_{12}, x_{21},x_{22}$ will be called, according to 
\cite{fj}, the  $quantum$ $Minkowski$ $space-time$. One can define the quantum Laplacian operator on this space by the formula:
\begin{eqnarray}
 \square_x=\p_{11}\p_{22}-\p_{12}\p_{21},
 \end{eqnarray}
 where $\p_{11}\p_{22},\p_{12},\p_{21}$ obey the following commutation relations:
 \begin{eqnarray}
 &&\p_{11}\p_{21} = \p_{21}\p_{11} ,\quad \p_{12}\p_{22} = \p_{22}\p_{21} , \nonumber\\
&&[\p_{11} , \p_{22} ] + [\p_{21} , \p_{12} ] = 0 , \nonumber\\
&&\p_{11}\p_{12} = q^{-2}\p_{12}x_{11} , \quad \p_{21}\p_{22} = q^{-2}\p_{22}\p_{21},
 \nonumber\\
&&\p_{12}\p_{21} = q^2\p_{21}\p_{12}.
 \end{eqnarray}
It was shown in \cite{fj} that the kernel of quantum Laplacian operator is spanned by elements of the following kind:
\begin{eqnarray}
X^{\lambda}_
{rk} =\frac{1}{2\pi i}
\oint (x_{11}s + x_{21})^{\frac{\lambda-r}{2}}(x_{12}s + x_{22} )^{\frac{\lambda+r}{2}}s^{\frac{k-\lambda}{2}-1} ds,
\end{eqnarray}
where 
\begin{eqnarray}
-\lambda \le r,k \le \lambda , \lambda\in{\mathbb{Z}_+},\quad r,k \equiv \lambda(\rm{mod} 2).
\end{eqnarray}
Moreover, the operator
\begin{eqnarray}
\tilde \square_x\equiv det(x)\square_x, \quad{\rm where} \quad det(x)=\delta^2det_q=x_{11}x_{22}-x_{12}x_{21}
\end{eqnarray}
 is diagonalizable on the quantum Minkowski space-time. 
 \begin{prop}{\rm\cite{fj}}
 The eigenvectors of the operator $\tilde \square_x$ have the form
\begin{eqnarray}
det(x)^jX^{\lambda}_
{rk}, 
\end{eqnarray}
such that the eigenvalues are 
\begin{eqnarray}\label{eig}
(j)_{q}(j+\lambda+1)_{q},
\end{eqnarray}
where $(n)_q=\frac{q^{2n}-1}{q^2-1}$.
 \end{prop}
\noindent Let us denote 
\begin{eqnarray}
Y^{\lambda,j}_{r,k}(A,B,C,D)\equiv \delta^{-2j-\lambda}det(x)^jX^{\lambda}_
{rk}.  
\end{eqnarray}
The following Proposition holds.
\begin{prop}
The elements $Y^{\lambda,j}_{r,k}(A,B,C,D)$ form a basis in the subspace $V_{\lambda,\lambda+2j}\otimes \bar{V}_{\lambda,\lambda+2j}$ of $GL^+_q(2)$.
\end{prop}
\noindent {\bf Proof.} It is enough to prove that $Y^{\lambda,0}_{r,k}(A,B,C,D)$ form a basis in $V_{\lambda,\lambda}\otimes \bar{V}_{\lambda,\lambda}$. 
We will prove this by induction in $\lambda$. The assertion is obvious in the case $\lambda=1$. Then $V_{1,1}\otimes \bar{V}_{1,1}$ has a basis $Y^{1,0}_{r,k}(A,B,C,D)$, where 
$r,k$ are equal to $0$ or $2$. Suppose that this statement holds for $\lambda-k$, where $\lambda\ge k\ge 0$. Let's prove it for $\lambda+1$. First of all, we know 
that  the monomials of $\lambda$ elements form a basis in the space  
$U_{\lambda}=\oplus_{\lambda\ge 2m\ge 0}V_{\lambda-2m,\lambda}\otimes \bar{V}_{\lambda-2m,\lambda}$. By induction, we know that $Y^{\lambda-2m}_{r,k}(A,B,C,D)$ form a basis in all $U_{\mu}$, where $\mu\le \lambda$.

As a consequence of Proposition 4.4. we have two natural actions of 
$U_q(gl(2))$ on $GL^+_q(2)$ inherited from vertex algebras. Let's consider the action of $U_q(gl(2))$, obtained from $\mathbb{G}_{-\varkappa,-\eta}$ braided VOA. The action of the 
corresponding $I, q^H,E,F$-generators on $A,B,C, D$ is as follows:
\begin{eqnarray}
&&IA=A\quad, IB=B,\quad IC=C,\quad ID=D,\nonumber\\
&&FB=A,\quad FD=C, \quad FC=0,\quad FA=0,\nonumber\\
&&q^HB=qB, \quad q^HD=qD,\quad q^HC=q^{-1}C,  \quad q^HA=q^{-1}A,\nonumber\\
&&EA=B,\quad EC=D, \quad EB=0, \quad ED=0.
\end{eqnarray}
One can see that action of $E,F$ interchanges the commuting elements 
$(x_{11}s + x_{21}) \leftrightarrow (x_{12}s + x_{22} )$ and therefore the space generated by $Y^{\lambda,0}_{r,k}(A,B,C,D)$ is invariant under $U_q(gl(2))$ action. Moreover, for fixed $k$, $Y^{\lambda,0}_{r,k}(A,B,C,D)$ span the $\lambda+1$
dimensional representation of $U_q(gl(2))$. At the same time, the elements   
of $U_{\lambda}$ which belong to $U_q(gl(2))$ irreducible representations 
with highest weight $\lambda$, should belong to the space $V_{\lambda,\lambda}\otimes \bar{V}_{\lambda,\lambda}$. Therefore, $Y^{\lambda,0}_{r,k}(A,B,C,D)\in 
V_{\lambda,\lambda}\otimes \bar{V}_{\lambda,\lambda}$ and hence $Y^{\lambda,0}_{r,k}(A,B,C,D)$ form a basis in $V_{\lambda,\lambda}\otimes \bar{V}_{\lambda,\lambda}$, since they are linearly independent. Thus the Proposition is proven.
\hfill$\blacksquare$\bigskip

\noindent Now we give meaning to the operator $\tilde\square_x$ in vertex algebra setting. 
Consider an operator $\tilde\square^{CFT}_x$ on the braided VOA $\hat{ \mathbb{G}}^+_{\varkappa,\eta}$ of the following form:
\begin{eqnarray}
\tilde\square^{CFT}_x=(\hat{a}_0-\alpha_0)_q(\hat{a}_0+1)_q,
\end{eqnarray}
where $\alpha_0$ is zero-mode of the quantum field $\alpha(z)$ from Feigin-Fuks realization of Virasoro algebra with generators $\{L_n\}$ and 
$\hat{a}_0=\eta^{-1}a_0$. 

Being constructed by means of zero-modes, this operator commutes with the action of $W_2$ algebra and semi-infinite cohomology operator. Comparing the eigenvalues of $\tilde\square^{CFT}_x$ and 
$\tilde\square_x$ (see (\ref{eig})) we obtain the following theorem.
\begin{theorem}
The operator $\tilde\square^{CFT}_x$ induces the operator $\tilde\square_x$ on \\
$(H^{\frac{\infty}{2}+0}(W_2,\mathbb{C}\mathbf{c},\hat{\mathbb{G}}^+),\mu)\cong
GL^+_q(2)$.
\end{theorem}

\section*{Acknowledgements}
We are indebted to P.I. Etingof, M. Jardim, A.A. Kirillov Jr., K. Styrkas, G.J. Zuckerman for fruitful discussions. 
The research of I.B.F. was supported by NSF grant DMS-0457444. A.M.Z. would like to thank the organizers of the Simons Workshop 2011, where this work was partly done.


\begin{thebibliography}{10} 
\bibitem{FFr} B. Feigin, E. Frenkel, {\it Integrals of Motion and Quantum Groups}, Lect. Notes in Math., 1620, Springer, Berlin 1996. 
arXiv:hep-th/9310022.
\bibitem{fb} E. Frenkel, D. Ben-Zvi, {\it Vertex algebras and algebraic curves}, AMS, Providence, USA (2004).
\bibitem{fgz} I.B. Frenkel, H. Garland, G.J. Zuckerman, {\it Semi-infinite cohomology and string theory}, Proc. Nat. Acad. Sci. {\bf 83} (1986) 8442-8446.
\bibitem{fj} I.B. Frenkel, M. Jardim, {\it Complex ADHM equations, sheaves on $\mathbb{P}^3$ and quantum instantons},  arXiv:math/0408027.
\bibitem{FHL} I.B. Frenkel, Y.-Z. Huang, J. Lepowsky, {\it On axiomatic approach  to vertex operator algebras and modules}, Memoirs of AMS {\bf 494} (1993).
\bibitem{fs}I.B. Frenkel, K. Styrkas, {\it Modified regular representations of affine and Virasoro algebras, VOA structure and semi-infinite cohomology}, Adv. Math. 206 (2006) 57-111.
\bibitem{fz}I.B. Frenkel, A.M. Zeitlin, {\it Quantum group as semi-infinite cohomology}, Comm. Math. Phys. {\bf 297} (2010) 687-732, arXiv:0812.1620.
\bibitem{lz1}B. H. Lian, G.J. Zuckerman, {\it 2D gravity with c=1 matter}, Phys. Lett. {\bf B266} (1991) 21-28.
\bibitem{lz1a}B. H. Lian, G.J. Zuckerman, {\it Semi-infinite cohomology and 2D Gravity. I}, Commun. Math. Phys. {\bf 145} (1992) 561-593. 
\bibitem{lz2}B. H. Lian, G.J. Zuckerman, {\it New Perspectives on the BRST-algebraic structure of string theory}, Commun. Math. Phys. {\bf 154} (1993) 613-646.
\bibitem{ms} G. Moore, N. Seiberg, {\it Classical and quantum conformal field theory}, Commun. Math. Phys. {\bf 123} (1989) 177-254.
\bibitem{styrkas} K. Styrkas, {\it Quantum groups, conformal field theories, and duality of tensor categories}, PhD Thesis, Yale, 1998. 
\end{thebibliography}
\end{document}